    \documentclass{amsart}
    \usepackage{amsfonts}
    \usepackage{amssymb}
    \usepackage{amsbsy}
\usepackage{amsthm}
\usepackage{amsmath}

\newcount\skewfactor
\def\mathunderaccent#1#2 {\let\theaccent#1\skewfactor#2
\mathpalette\putaccentunder}
\def\putaccentunder#1#2{\oalign{$#1#2$\crcr\hidewidth
\vbox to.2ex{\hbox{$#1\skew\skewfactor\theaccent{}$}\vss}\hidewidth}}
\def\name{\mathunderaccent\widetilde-3 }

\makeatletter
\makeatletter

\newcommand{\type}{\operatorname{type}\nolimits}
\newcommand{\supp}{\operatorname{supp}}
\newcommand{\dom}{\operatorname{dom}}
\newcommand{\ran}{\operatorname{ran}}
\newcommand{\stronger}{\ge}
\newcommand{\weaker}{\le}
\newcommand{\liff}{\Leftrightarrow}
\newcommand{\limpl}{\Rightarrow}

\newcommand{\F}{{\mathcal F}}

\newcommand{\M}{{\mathfrak M}}
\newcommand{\X}{{\mathfrak X}}
\newcommand{\N}{{\mathfrak N}}
\newcommand{\p}{{\mathfrak p}}
\newcommand{\q}{{\mathfrak q}}
\renewcommand{\c}{{\mathfrak c}}
\renewcommand{\r}{{\mathfrak r}}
\renewcommand{\P}{{\mathfrak P}}
\newcommand{\Q}{{\mathfrak Q}}
\newcommand{\QQ}{{\mathbb Q}}
\newcommand{\PP}{{\mathbb P}}

\newcommand{\<}{\sqsubseteq}

\newcommand{\forces}{\Vdash}

\newcommand{\amal}[2]{ \oplus_{#1, #2} } 

\newcommand{\claim}[1]{\par\noindent{\bf #1.} }

\newcommand{\on}{{\upharpoonright}}

\newenvironment{ITEMIZE}%
        {\begin{list}{$*$}{\labelwidth1.6cm\leftmargin2cm\rightmargin1cm
                \itemsep0.1cm\relax\let\item\itemizeitem}}%
        {\end{list}}

\makeatletter
\def\enumerate{\ifnum \@enumdepth >3 \@toodeep\else
      \advance\@enumdepth \@ne
      \edef\@enumctr{enum\romannumeral\the\@enumdepth}\list
      {\csname label\@enumctr\endcsname}{\parsep0pt\itemsep0pt
\usecounter
        {\@enumctr}\def\makelabel##1{\hss\llap{##1}}}\fi}

\def\itemize{\ifnum \@itemdepth >3 \@toodeep\else \advance\@itemdepth \@ne
\edef\@itemitem{labelitem\romannumeral\the\@itemdepth}%
\list{\csname\@itemitem\endcsname}{\parsep0pt\itemsep0pt
\def\makelabel##1{\hss\llap{##1}}}\fi}

\renewcommand{\cases}[1]{
\left \{\,\vcenter
 {\normalbaselines \m@th \ialign {$##\hfil $&\quad ##\hfil
 \crcr #1\crcr }}\right .}
\makeatother

\swapnumbers

\def\xx#1 {%
\newtheorem{#1}[thm]{#1}}
\def\yy#1 {%
\newenvironment{#1}{\begin{roman-#1}\rm}{\end{roman-#1}}
\newtheorem{roman-#1}[thm]{#1}}

\yy Fact
\xx Theorem
\xx Corollary
\xx Lemma
\yy Definition
\yy Remark
\yy Notation
\yy {Fact and Definition}
\yy Acknowledgments
\yy Setup
\yy Note
\yy {Overview of the construction}
\xx Claim
\xx Conclusion

\def\h0{the strong chain condition}
\def\m0{creature}

\begin{document}
\title{A Partial Order Where All Monotone Maps Are Definable}
\author{Martin Goldstern \and Saharon Shelah}
\thanks{This is number 554 in the second author's
        publication list. Supported by the Israeli Academy of Sciences.}
\address{Technische Universit\"at\\Wiedner Hauptstra\ss{}e
  8--10/118.2\\
A-1040 Wien}
\curraddr{Mathematik WE2 \\Freie Universit\"at\\
Arnimallee 3\\D-14195 Berlin}
\email{goldstern@tuwien.ac.at}
\address{Department of Mathematics\\Hebrew University, Givat Ram\\
91904 Jerusalem, Israel}
\email{shelah@math.huji.ac.il}

\date{February 1997}

\begin{abstract}
It is consistent that there is a partial order $(P,{\le})$ of size 
$ \aleph_1  $
such that every monotone function $f:P\to P$ is first order definable
in $(P,{\le})$
\end{abstract}

\maketitle

\makeatletter
\let\old@item\@item
\def\itemizeitem{\@ifnextchar[{\itemize@item}{\itemize@item[\@itemlabel]}}
\def\itemize@item[#1]{\old@item[{#1}]\if@noparitem\else
                        \def\@currentlabel{#1}\fi}
\makeatother

It is an open problem whether there can be an infinite lattice $L$
such that every monotone function from $L$ to $L$ is a polynomial. 
Kaiser and Sauer \cite{KS} showed that such a lattice would have to be
bounded, and cannot be countable. 

Sauer then asked the weaker question if there can be an
infinite partial order $(P,\le)$ such that all monotone maps from $P$ to $P$
are at least definable. (Throughout the paper, 
``definable'' means ``definable with parameters by a
first order formula in the structure $(P,\le)$.)

Since  every infinite partial order $P$ 
 admits  at least 
$\c=2^{\aleph_0}$ many monotone maps from $P$ to $P$, our partial 
order must have size (at least) continuum.

We show:

\begin{Theorem}  The statement 
\begin{quote}
There is a partial order $(P, {\le})$ of size $\c={\omega_1} $ such that all
monotone functions $f:P \to P$ are definable in $P$
\end{quote}
is consistent relative to ZFC.  Moreover, the statement holds in any
model obtained by adding (iteratively) ${\omega_1}$ Cohen reals to a
model of CH. 
\end{Theorem}

We do not know if Sauer's question can be answered outright (i.e., in
ZFC), or even from CH.

\subsection*{Structure of the paper}


In section \ref{four} we  give four conditions
on a partial order on $(P, {\<})$ of size $\kappa$
  and we show that they are
sufficient to ensure the conclusion of the theorem.  This section is
very elementary.  The two main conditions of section \ref{four} are
\begin{enumerate}
\item a requirement on small sets, namely that they should be definable
\item two requirements on large sets (among them: ``there are no large
        antichains'')
\end{enumerate}
Here, ``small''
means of size $< |P|$, and ``large'' means of  size $= |P|$. 


 In section \ref{code} we show how to take care
        of requirement 1 in an inductive construction of our partial
order in $\c$ many steps. Each definability requirement will
be satisfied at some  stage $\alpha<\c$.

Finally, in section \ref{forcing} we deal with the problem of avoiding large
antichains.   Here the 
 inductive construction is not so  straightforward, as we
have to ``anticipate'' potential large sets and ensure that in the end
they will not contradict our requirement.    The standard tool for
dealing with such a problem is $\diamondsuit$.   This combinatorial
principle has been used for a related construction in \cite{Sh:128},
and it is possible to use the techniques of \cite[section 5]{Sh:136}
to combine it with the requirement on small sets to show that the
conclusion of our theorem actually follows from $\diamondsuit$.
  However, we use instead a forcing construction which seems to
  be somewhat simpler. 
 We will work in an iterated forcing extension, and
use a ${\Delta} $-system argument to ensure that the requirements
about 
large sets are met.    This argument is carried out in section
\ref{forcing}, 
which therefore requires a basic knowledge of forcing.

\section{Four conditions}
\label{four}

\begin{Theorem}\label{4prop}
Assume that $(P, {\<})$ is a partial order, $\kappa:=|P|$ a regular cardinal. 
We will call subsets of cardinality $<\kappa$ ``small.''
Assume that  the following conditions hold: 
\begin{ITEMIZE}
\item[C1] Every antichain is small
(an antichain is a set of
        pairwise incomparable elements).
\label{ccc}
\item[C2] Whenever $g:(P, {\<}) \to (P,{\<})$ is monotone, then
        there is a small  set $A \subseteq P$ such that for
        all $\alpha\in P$, $f(\alpha)\in A \cup \{\alpha\}$.
\label{bounded}
\item[C3] 
For all $\alpha \in P$ the set 
        $\{\beta\in P: \beta \< \alpha\}$ 
        is small. 
\label{omega1like}
\item [C4] Every small  subset of $P \times P$ is
        definable in the structure 
        $(P, {\<}) $. 
\label{defble}
\end{ITEMIZE}
Then every monotone map from 
        $(P, {\<}) $ to itself is  definable in
        $(P, {\<}) $.  
\end{Theorem}

We will prove this theorem below.

\begin{Lemma}\label{aprime}
Let $\<$ be as above.   Then for any set $A  \subseteq P$ 
there is a small  set $A' \subseteq A$ 
such that 
 $\forall \alpha \in A \, \exists {\gamma} \in A':  {\gamma} \<
        \alpha$. 
\end{Lemma}
\begin{proof}
  Let $ B \subseteq A$ be a maximal antichain in $A$.  So $B$ 
is small, by \ref{ccc}.  Let $A'$ be the downward closure 
of $B$ in $A$, i.e., 
$A':= \bigcup\limits_{ \beta  \in  B} 
        \{\gamma\in A: \gamma \< \beta\}$.
Then $A'$ is still small because of \ref{omega1like}.  
Clearly $A'$ is as required. 
\end{proof}

\begin{Fact}\label{2.3}
Let $(P, {\<})$ be any partial order, $g:P \to P$ a monotone 
function.  Then $g$ is definable iff the set 
$\{ (x,y): g(x) \< y \}$ is definable. 
\end{Fact}
\begin{proof}  Let $B:= \{ (x,y): g(x) \< y \}$.  Clearly $B$ is definable
from $g$.   Conversely, $g$ can be defined from $B$ as follows: 
$$g(x) = z  \ \liff \ \left(\forall y \, (x,y) \in B \limpl z \<
y\right)
\text { and } 
\left((x,z)\in B\right) $$
\end{proof}

\subsection*{Proof of \ref{4prop}:} Let $g : P \to P$ be
monotone. Let $P_0 \subseteq P $ be small   
such that \begin{itemize}
\item 
        If $g(\alpha) \not= \alpha$ then  $g(\alpha)\in P_0$. 
\item  $\beta\< \alpha \in P_0 $ implies 
        $\beta\in P_0$. 
\end{itemize}
(Such a set $P_0$ can be found using \ref{bounded} and
\ref{omega1like}). 
For every $j \in P_0$ let $D_j:=
\{ \alpha \in P \setminus P_0 : g(\alpha) = j\}$, and let
$D:= \{\alpha \in P\setminus P_0: g(\alpha)= \alpha\}$. 
We can find  $D'$ and $D_j'$ as in \ref{aprime}. 
Let $P_1 = \bigcup_j D'_j \cup D'$.

\claim{Claim 1} For $\alpha \in P\setminus P_0$ we have:
$$ g(\alpha) = \alpha \ \Longleftrightarrow \ \exists {\beta} \in P_1 : 
{\beta} \< \alpha,  g(\beta)= \beta 
$$ 
The direction $\Rightarrow$ follows immediately from the definition of
$D'$. For the other direction, assume that ${\beta} \in P_1$, ${\beta}
\< \alpha$, $g(\beta)=\beta$.   Since $g(\beta) \< g(\alpha)$, 
and $g(\beta)=\beta\notin P_0$, we also have $g(\alpha)\notin P_0$.  Hence  
$g(\alpha)=\alpha$. 

\claim{Claim 2}
For   $\alpha \in P \setminus P_0$, $g(\alpha) \not= \alpha $ we have: 
$$ g(\alpha ) \< i \ \Longleftrightarrow \ 
\forall {\gamma}\in P_1: {\gamma} \< \alpha \limpl  g(\gamma) \< i$$
The implication $\Rightarrow $ follows from the monotonicity of $g$. 
For the converse direction, assume that $g(\alpha)=j$.  
Since $D'_j \subseteq P_1$, we can find $\gamma\in P_1$, $\gamma\<\alpha$ 
and $g(\gamma) = j$. By assumption, $g(\gamma)  \< i$, 
so $g(\alpha) \< i$.

\medskip
Claims 1 and 2, together with \ref{2.3} now imply that 
 $g$ can be defined from the graphs of the functions  $g\on P_0 $ and
 $g \on P_1$.

\section{Coding small sets}\label{code}

In this section we will show how to build a partial order satisfying
C3 and C4.   Throughout this section we assume the continuum
hypothesis.

\begin{Definition}
\begin{enumerate}
\item
We call $\M$ a ``\m0'', if 
$\M = (M,{\<}, F,H)$, where
\begin{itemize}
\item $\<$ is a partial order on $M$
\item $F$ is a partial symmetric function, $\dom(F) \subseteq M \times
        M$, $\ran(F) \subseteq M$
\item if $F(x,y)=z$, then  $x  <  z $, $y< z $, and there is no $z'<z$
with $x\le z'$ and $y\le z'$ (i.e., $z$ is a minimal upper bound of
$x,y$, and $x,y$ are incomparable).   
\item $F$ is locally finite, i.e.:  For any finite $A \subseteq M$
there is a finite set $B$, $A \subseteq B \subseteq M$ such that for
any $(x,y)\in \dom F \cap (B \times B)$, $F(x,y)\in B$. 
\item $H \subseteq M\times M \times M$.   We define 
        $H(x,y):=\{z: H(x,y,z)\}$
\item If $H(x,y,z)$ then $z$ is a minimal upper bound of $x,y$, and
$x,y$ are incomparable
\item  $H(x,y)=H(y,x)$
\item If $H(x,y,z) $ then $(x,y)\notin \dom(F)$. 
\end{itemize} 
\item If $\M_1$, $\M_2$ are \m0s, then $\M_1 \le \M_2$ ($\M_2$ is an
extension of $\M_1$) means that  
$M_1 \subseteq M_2$, and the relations/functions of $\M_1$ are just
the restrictions of the corresponding relations/functions of $\M_2$. 
(so in particular, $\dom(F^{\M_2}) \cap (M_1 \times M_2) =
\dom(F^{\M_1})$. 
\item 
We say that $\M_2$ is an ``end extension'' of $\M_1$ (or that $\M_1$
is an ``initial segment'' of $\M_2$) if $\M_1 \le \M_2$ and  
$\M_1$ is downward closed in
$\M_2$, i.e., $\M_2 \models x\< y$, $y \in M_1$ implies $x\in M_1$. 
\item We use the above terminology also if one or both of the
structures $\M_1$, $ \M_2$ is just a partial order without an
additional function $F$ or relation $H$. 
\end{enumerate}
\end{Definition}

\begin{Notation}
When we consider several \m0s $\M$, $\M_1$, $\M^*$, \dots,   then it is
understood that their universes are called $M$, $M_1$, $M^*$, \dots,
their partial orders are $\<$, $\<_1$, $\<^*$, \dots, etc. 
Instead of writing $x \<_1 y$ we may write $\M_1 \models x \< y$, etc. 

We will use the letters $\M$ and $\N$ to denote possibly infinite
\m0s, and $\p$, $\q$ to denote finite \m0s. 

\end{Notation}

\begin{Definition}\label{fdef}
If $(M, {\<})$ is a partial order, then we let $\F({\<})$ be the partial
binary function $F$ satisfying 
$$ F(x,y) = 
\left\{ \ 
\parbox[c]{8cm}{the unique minimal upper bound of $x$ and $y$,
if it exists and if $x$ and $y$ are incomparable \\[3mm]
undefined, otherwise }
\right.$$
\end{Definition}

\begin{Setup}\label{setup}
Let $(M_\delta: \delta<\c)$ be a continuous increasing sequence of
infinite sets    with $|M_{\delta+1} \setminus M_\delta| =
|M_\delta| < \c$.
Let $(R_\delta: \delta < \c)$ be a sequence of relations
with $R_\delta \subseteq M_\delta \times M_\delta$ such that for any
$R \subseteq M_\delta \times M_\delta$ there is a $\delta ' > \delta$
such that $R= R_{\delta '}$.    (Such a sequence can be found since we
are assuming CH).

For each $\delta$ we fix some element $e_\delta\in M_{\delta+1}
\setminus 
M_\delta$. For each $(\alpha,\beta)\in M_\delta\times M_\delta$
we pick disjoint sets 
$A_{\alpha\beta}^\delta  $, $B_{\alpha\beta}^\delta $,
$C_{\alpha\beta}^\delta $, 
$\Delta_{\alpha\beta}^\delta $, $\Gamma_{\alpha\beta}^\delta $
(we will omit the superscript $\delta$ if it is clear from the context)
 satisfying
\begin{gather}
  |A_{\alpha\beta} | = 2, \ |B_{\alpha\beta} | = 3,\  |C_{\alpha\beta} | = 1\\
  |\Delta_{\alpha\beta} | = 3,  \ 
         \Delta_{\alpha\beta}=\{a_{\alpha\beta},
                b_{\alpha\beta},c_{\alpha\beta}\},  \ 
  \Gamma_{\alpha\beta}  =
                \{\gamma_{\alpha\beta}\} 
\end{gather}
where $e_{\delta} $ is not in any of these sets and the sets 
$$ \Omega _{\alpha\beta}  := 
A_{\alpha\beta} \cup B_{\alpha\beta} \cup C_{\alpha\beta} \cup
        \Delta_{\alpha\beta} \cup \Gamma_{\alpha\beta}  $$
satisfy $(\alpha,\beta)\not= (\alpha ', \beta ') \limpl
\Omega_{\alpha\beta} \cap \Omega _{\alpha'\beta'} =\emptyset$. 
For $x\in \Omega_{\alpha\beta} $ we define $\Omega(x):= 
\Omega_{\alpha\beta} \cup \{\alpha,\beta\}$. 

Moreover, we choose the sets $\Omega _{\alpha\beta} $ such that
$M_{\delta+1} \setminus  \bigcup_{\alpha,\beta\in M_\delta}
\Omega_{\alpha\beta}$ 
still has size $|M_\delta|$. 
For $ x \in M_{\delta+1} \setminus  \bigcup_{\alpha,\beta \in
M_\delta} 
\Omega_{\alpha\beta}$ we let $\Omega(x)=\emptyset$. 
\end{Setup}

\begin{Definition}
\begin{enumerate}
\item Let $\M$ be a \m0. We say that $\Delta \in [M]^3$ is a
{\em triangle} (also called an $\M$-triangle or $F$-triangle if
$F=F^\M$) iff $\dom(F^\M) \supseteq [\Delta]^2$. 
\item Let $\Delta$ be an $F$-triangle.  We say that $a\in
\Delta $ is a {\em base point} for $\Delta$ if there is a
unique (unordered) pair $\{b,c\}$ with $F(b,c)=a$.  
\item Let $\Delta$ be a triangle.  We say that $b$ is
{\it anchor} for $\Delta$ if there is a $c$ such that $F(b,c)$
is a base point for $\Delta$. 
\item If $\M_1 \le \M_2 $ then we say that $\M_ 2 $ is a ``separated''
extension of $\M_1$ iff 
\begin{itemize}
\item[$\circ$] every triangle in $\M_2$ is contained in $M_1$
        or in $M_ 2 \setminus M_1$, and
\item[$\circ$] if $\Delta  \subseteq M_ 2 \setminus M_1$ is a
        triangle, then $\Delta$ is not anchored at any point in
        $M_1$.
\item[$\circ$] If $\M_2 \models H(x,y,z)$, and $x,y\in M_1$ then also
        $z\in M_1$. 
\end{itemize}
\end{enumerate}
\end{Definition}

\begin{Note}
  Our goal is to
define $F_{\delta+1} $ such that $R_\delta$ will become definable
(and at the same time make it possible for $F_{\delta+1}$ to be of the
form $\F({\<})$ for some end extension $\<$ of $\<_\delta$). We will
achieve this by ``attaching'' triangles to all pairs $(\alpha, \beta)$
in $R_\delta$ in a way that 
 $(\alpha,\beta)$ can be reconstructed from the
triangle.   We also need to ensure that the only triangles in
        $\M_{\delta+1} $ (and also in any $\M_{\delta '}$, $\delta ' >
        \delta$) are the ones we explicitly put there. 

The particular way of coding pairs by triangles is rather arbitrary. 
\end{Note}

\begin{Overview of the construction}\label{construct}
By induction on $\delta\le \c$ we will define \m0s $\M_\delta = (M_\delta,
\<_\delta, F_\delta, H_{\delta} )$ such that 
\begin{ITEMIZE}
\item[(A)] $\gamma < \delta $ implies that $\M_\delta$ is
        a separated end extension of $\M_\gamma$
\item[(B)] $F_\delta = \F(\<_\delta)$.  (See \ref{fdef})
\item[(C)] $R_\delta$ is definable in $(M_{\delta+1}, {\<_{\delta+1}},
        F_{\delta+1} )$
        and hence also in $(M_{\delta+1},{\<_{\delta+1}})$, and the
        definition of $R_\delta$ is absolute for any separated end
        extension of $\M_{\delta+1} $. \end{ITEMIZE}
For limit $\delta $ we let $\<_\delta = \bigcup_{\gamma < \delta}
\<_\gamma$, $F_\delta= \bigcup_{\gamma < \delta } F_\gamma$, 
 $H_\delta= \bigcup_{\gamma < \delta } H_\gamma$. 
\end{Overview of the construction}

We will construct $\M_{\delta+1} $ from $\M_\delta$ in two
steps.  First we define a  function $F_{\delta+1}$ 
such that $R_\delta$ becomes definable in $(M_{\delta+1},
F_{\delta+1})$.  Then we show that we can find a partial order
$\<_{\delta+1}$ (end-extending $\<_\delta$) such that 
$F_{\delta+1}=\F({\<_{\delta+1}})$.

\subsection*{Construction of $F_{\delta+1} $}
$F_{\delta+1}$ will be defined as follows:
\begin{ITEMIZE}
\item $F_{\delta+1} \on (M_\delta\times M_\delta) = F_\delta $
\item If $(\alpha,\beta)\in R_\delta$ then
\begin{ITEMIZE}
        \item $F(\alpha,x)= a_{\alpha\beta}$ for $x\in
                A_{\alpha\beta}$
        \item $F(\beta,x)= b_{\alpha\beta}$ for $x\in
                B_{\alpha\beta}$
        \item $F(e_\delta ,x)= c_{\alpha\beta}$ for $x\in
                C_{\alpha\beta}$
        \item $F(x,y)=\gamma_{\alpha\beta}$ for $\{x,y\}\in
                [\Delta_{\alpha\beta}]^2 $. 
\end{ITEMIZE}
\item Except where required by the above  (and by symmetry),
                $F_{\delta+1} $ is undefined.
\end{ITEMIZE}

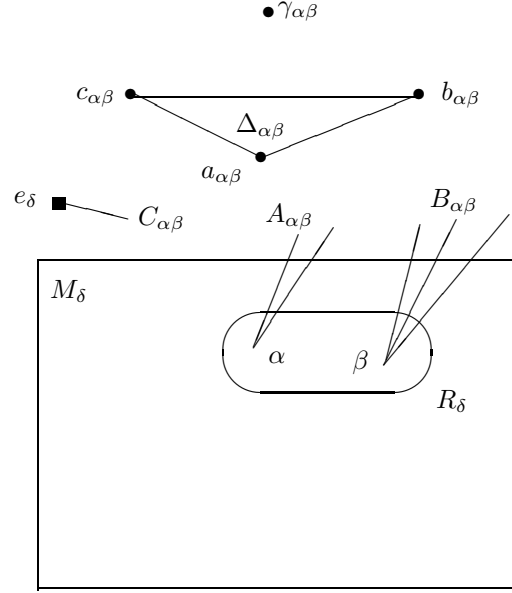
\unitlength 1.00mm
\linethickness{0.4pt}
\begin{picture}(105.33,116.41)
\put(41.33,38.67){\line(0,1){44.00}}
\put(41.33,82.67){\line(1,0){64.00}}
\put(105.33,82.67){\line(0,-1){43.67}}
\put(105.33,39.00){\line(-1,0){64.00}}
\put(53.66,104.66){\circle*{1.49}}
\put(71.00,96.33){\circle*{1.49}}
\put(71.00,96.33){\line(-2,1){17.33}}
\put(53.66,104.33){\line(1,0){38.67}}
\put(71.33,96.33){\line(5,2){21.00}}
\put(92.00,104.66){\circle*{1.49}}
\put(41.33,90.67){\makebox(0,0)[rc]{$e_\delta$}}
\put(44.33,90.33){\line(4,-1){9.00}}
\put(54.67,88.00){\makebox(0,0)[lc]{$C_{\alpha \beta}$}}
\put(76.00,86.00){\line(-2,-5){6.00}}
\put(70.00,71.00){\line(2,3){10.67}}
\put(104.00,88.67){\line(-5,-6){16.67}}
\put(87.33,68.67){\line(1,2){9.67}}
\put(87.67,69.33){\line(1,4){4.50}}
\put(79.83,70.33){\oval(27.67,10.67)[]}
\put(94.33,65.33){\makebox(0,0)[lt]{$R_\delta$}}
\put(72.00,69.67){\makebox(0,0)[lc]{$\alpha$}}
\put(85.33,69.00){\makebox(0,0)[rc]{$\beta$}}
\put(71.00,100.33){\makebox(0,0)[cc]{$\Delta_{\alpha\beta}$}}
\put(72.00,115.66){\circle*{1.49}}
\put(76.00,115.67){\makebox(0,0)[cc]{$\gamma_{\alpha\beta}$}}
\put(77.99,88.34){\makebox(0,0)[rc]{$A_{\alpha\beta}$}}
\put(96.66,90.34){\makebox(0,0)[cc]{$B_{\alpha\beta}$}}
\put(43.00,80.00){\makebox(0,0)[lt]{$M_\delta$}}
\put(51.66,104.33){\makebox(0,0)[rc]{$c_{\alpha\beta}$}}
\put(68.66,95.00){\makebox(0,0)[rt]{$a_{\alpha\beta}$}}
\put(95.00,104.66){\makebox(0,0)[lc]{$b_{\alpha\beta}$}}
\put(43.33,89.33){\rule{1.67\unitlength}{1.67\unitlength}}
\end{picture}
\vskip -3cm

This completes the definition of $F_{\delta+1}$.     The diagram above
is supposed to illustrate this definition.   Pairs $(x,y)$ on which
$F_{\delta+1}$  is defined are connected by a line; the value of
$F_{\delta+1}$ at such pairs  is the small black disk  above the pair.

\begin{Fact}\label{cfact}
\begin{enumerate}
\item \label{def0}
$(M_{\delta+1}, F_{\delta+1}) $ is a separated extension of
$\M_\delta$. 
\item \label{def1}
$R_\delta$ is definable in $(M_{\delta+1} , F_{\delta+1} )$. 
\item \label{def2}
Let $(M_{\delta+1}, F_{\delta+1}) \le (M,F)$ and assume that
\begin{ITEMIZE}
\item $(M,F)$ is a separated extension of $(M_{\delta+1} ,
F_{\delta+1})$
%
\item If $F(x,y) \in M_{\delta+1}$ then $x,y\in M_{\delta+1}$. 
(This is certainly true if $F=\F(\<)$, where $\<$ is an end extension
of $\<_{\delta+1} $.)
\end{ITEMIZE}
Then $R_\delta$ is definable in $(M,F)$. 
\end{enumerate}
\end{Fact}

\begin{proof}  \ref{def0} is clear. 
 \ref{def1} is a special case of \ref{def2}.   Let $(M,F)$ be
as in \ref{def2}. 
Then      $R_\delta $ is the set of all pairs $(\alpha,\beta)$ such that
                 there is a triangle $\Delta$ 
                with a unique base, anchored at $e_\delta$ such that 
\begin{enumerate}
\item[either] $\alpha \not=\beta $ and
        $|\{x:F(\alpha,x)\in \Delta \}|=2$, 
        $|\{x:F(\beta ,x)\in \Delta \}|=3 $. 
\item[or] 
$\alpha =\beta $ and
        $|\{x:F(\alpha,x)\in \Delta \}|=5$,
\end{enumerate}
because the only triangles anchored at $e_\delta$ will be the sets
$\Delta_{\alpha\beta}^\delta $.   Clearly this is a definition in
first order logic with the parameter $e_\delta$. 
\end{proof}

\subsection*{Construction of $\<_{\delta+1} $} 
We will define $\<_{\delta+1} $ from a (sufficiently) generic filter for a
forcing notion $\QQ_\delta =  \QQ(\M_\delta, F_{\delta+1})$.

\begin{Definition}
Assume that $\M_\delta$, $F_{\delta+1} $ are as above.  We define the
forcing notion $\QQ_\delta$ as the set of all $\p$ such that
\begin{ITEMIZE}
\item[--] $\p$ is a finite \m0, $\p= (p, \<_\p, F_\p, H_\p)$
\item[--] $\p \on M_\delta \le \M_\delta$
\item[--] $p \subseteq M_{\delta+1}$,
\item[--] For all $x \in p \setminus M_{\delta}$, $\Omega(x) \subseteq p$. 
\item[--] $F_\p = F_{\delta+1} \on (p\times p)$
\item[--] $\p$ is a separated end extension of $\p \on M_\delta$.  
\end{ITEMIZE}
\end{Definition}

Clearly $(\QQ_\delta,{ \le})$ is a partial order. 
{\bf Note:  We force ``upwards,'' i.e., $\p \le \q$ means that $\q$ is
a stronger condition than $\p$.}  But we still call the generic set
``filter'' and not ``ideal.''

 If $G \subseteq Q_\delta$ is
a filter then there is a smallest \m0 which extends all $\p\in G$.
We call this \m0 $\M_G$.

\begin{Fact}\label{dense1}
The following sets are dense in $\QQ_\delta$. 
\begin{ITEMIZE}
\item[(a)] $D_x:=\{\p: x \in p \}$, for any $x\in M_{\delta+1} $.
\item[(b)] $D_{\p,x,y}:=\{\q: \q \perp \p \hbox{ or  } 
\q \models  \exists z \notin p: H(x,y,z)\}$, whenever $x,y\in p$ are
incomparable and
$(x,y)\notin \dom(F_\p)$.
\item[(c)] $E_{\p,x,y,z}:=\{\q: \q \perp \p \hbox{ or  }
\q \models  \exists z'\< z:   x\< z', y \< z', z\not=z'\}$,
 whenever $x,y\in p$, $(x,y)\in \dom(F_\p)$, $x\<z$, $y\<z$, 
$z\not= F(x,y)$.
\end{ITEMIZE}
\end{Fact}

\begin{Fact}
If $G$ meets all the dense sets above, then $\M_{\delta+1} := \M_G $
 is a separated 
end extension of $\M_\delta$, $F_{\delta+1} = \F({\<_{\M_G}})$.  
\end{Fact}
\begin{proof} We only prove the last statement. Let $F = \F({\<_{\M_G}})$.  
First assume  $F_{{\delta}+1}(x,y) = z^*$.  Clearly $z^*$ is a minimal
upper bound of $x,y$, so to prove $F(x,y)=z^*$ 
we only have to show that $z^*$ is the {\it  unique} minimal upper bound.
If
$z\not=z^*$ is also a minimal upper bound then we can find $p\in G$
containing $\{x,y,z,z^*\}$.  Now use \ref{dense1}(c) to find
$z'\not=z$ in 
$M_{\delta+1}$, $x,y \<z' \< z$, which is a contradiction. 
\\   
Now assume that $F_{\delta+1}(x,y)$ is undefined.   We have to show that 
also $F(x,y)$ is undefined. 
Applying \ref{dense1}(b) twice we can find two distinct
elements $z_1$, $z_2$ such that $H(x,y,z_1)$ and $H(x,y,z_2)$ both
hold in $\M_{\delta+1}$, hence there is no unique minimal upper bound
of $x$ and $y$, so  $F(x,y)$ is undefined.   
\end{proof}

\begin{Fact}
$\QQ_\delta$ is  countable, hence  satisfies the ccc. 
\end{Fact}

\begin{Conclusion}\label{conclusion}
 CH implies that a sufficiently generic  filter exists. 
Thus we have completed the definition of $\<_{\delta+1}$ and 
$H_{{\delta} +1}$. 
Clearly $\M_{\delta+1} $ will be as required. 
\end{Conclusion}

  In the last 
section we will 
show how to embed this construction into an iterated forcing argument.

\section{Amalgamation} \label{amalgamation}

Starting from a model of GCH we will construct an iterated forcing
notion $(\PP_\delta, Q_\delta: \delta < \kappa)$ such that each  partial
order $Q_\delta$ will be some $\QQ(\M_\delta, F_{\delta+1} )$ as in
the previous section.  This will define us a model $\M_\kappa$.  An
additional argument is then needed to show that $\M_\kappa$ will
satisfy \ref{ccc} and \ref{bounded}.  In this section we prepare some
tools for this additional argument by collecting some facts about
amalgamation. 

\bigskip

\begin{Definition}
Let $\p$ and $\q$ be (finite) \m0s,  $x\in p$, $y \in q$. We define 
$\p \oplus \q$ and $\p \amal xy \q$ as follows: 
\begin{ITEMIZE}
\item[] $\p \oplus \q = 
  (p \cup q, 
  {\<_{\p\oplus \q}}, 
  F_{\p\oplus \q},
    H_{\p\oplus \q})$,
  where $\<_{p\oplus q}$ is the transitive closure of
        $\<_\p\cup \<_\q$, and $F_{\p\oplus\q} = F_\p \cup F_\q$,
        $H_{\p\oplus\q} = H_\p \cup H_\q$. 
\item[] $\p \amal{x}{y} \q = (p \cup q, \<_{\p\amal{x}{y} \q},
        F_{\p\amal{x}{y} \q},
            H_{\p\amal xy \q})$,
         where $\<_{\p\amal{x}{y}\q}$ is the transitive closure of
        $\<_\p\cup \<_\q \cup \{(x,y)\}$, and
        $F_{\p\amal{x}{y}\q} = F_\p \cup F_\q$.  
        $H_{\p\amal{x}{y}\q} = H_\p \cup H_\q$.  
\end{ITEMIZE}
\end{Definition}

\begin{Fact}\label{amalfact}
 \begin{enumerate}
\item Assume that $\<_\p $ and $\<_\q$ agree on $p \cap q$, and
        similarly $F_\p$ and $F_\q$.   Then $\p \oplus \q$ is a \m0,
        $\p \< \p \oplus \q$, $\q \< \p \oplus \q$. 
\item If $\p$ and $\q$ are as above, and moreover: $\p$ and $\q$ are
        separated end extensions of $\r:=\p \cap \q$, and 
$ \type_\p(x,\r) = \type_\q(y,\r)$ (that is, for any $z\in r$ we have
        $\p \models z\< x $ iff $\q \models z \< y$), then also $\p
        \amal xy \q $ is a model extending $\q$ and end-extending
        $\p$, and
$$
(*) \ \ 
\begin{array}{rcl}
\p \amal xy \q \models a \< b  & \liff &      \p \models a \< b \\
                &&      \text{or } \q \models a \< b \\
                &&      \text {or } \p \models a \< x, \ \q \models y \< b
\end{array}
$$
\end{enumerate}
\end{Fact}

%
%
%
%
%
%
%
%
%
%
\begin{proof} We leave (1) to the reader. 
Let $\<^*$ be the relation defined in $(*)$.  First we have to check
that $\<^*$ is transitive.  Note that the
third clause in the definition of $a\<^* b$ can only apply if $a\in p$
and $b \in q \setminus r = q \setminus p$. 
\\
Let $a \<^* b \<^* c$.   {\it A priori}, there are 9 possible cases:

\def\casep#1#2{ #1 \<_\p #2 }
\def\caseq#1#2{ #1 \<_\q #2 }
\def\casex#1#2{#1 \<_\p x,\ y \<_\q #2}
\begin{tabular}{llcl}
$\casep ab$     &$\casep bc$ & $\Rightarrow$ &
 $\casep  a  c$ \\
$\casep ab$     &$\caseq bc$ & $\Rightarrow$ &
 $b\in p \cap q$, $\caseq  a  b$, so $\caseq a  c$ \\
$\casep ab$     &$\casex bc$ & $\Rightarrow$ &
 $\casex a c$ \\
$\caseq ab$     &$\casep bc$ & $\Rightarrow$ &
 $b\in p \cap q$, $\casep  a  b$, so $\casep a  c$      \\
$\caseq ab$     &$\caseq bc$ & $\Rightarrow$ &
 $\caseq  a  c $        \\
$\caseq ab$     &$\casex bc$ & $\Rightarrow$ &
 $b\in p \cap q$, $\casep  a b$, $\casex a c$                   \\
$\casex ab$     &$\casep bc$ & $\Rightarrow$ &
 impossible             \\
$\casex ab$     &$\caseq bc$ & $\Rightarrow$ &
 $\casex a c $                                                  \\
$\casex ab$     &$\casex bc$ & $\Rightarrow$ &
 impossible             \\
\end{tabular}

So we see that in any case we get $a \<^* c$.

Clearly $\p \amal xy \q $ is an end extension of $\p$.  We now check
that $ \q \le \p \amal xy \q$. Let $a,b\in q$, $ \p \amal xy \q
\models a \< b$.  The only nontrivial case is that $\p \models a \<
x$, $\q \models y \le b$.   Since $x $ and $y$ have the same type over
$\p \cap\q$, we also have $\q \models a \< y$.  Hence $\q \models a \<
b$, so $ \q \le \p \amal xy \q$. 

 Finally we check that $\p \amal xy \q$ is a \m0. 
 It is clear that $F_{\p \amal xy \q}$ is locally finite. \\
To check that the conditions on $H_{\p \amal xy \q}$ and $F_{\p \amal
  xy \q}$  (in particular, minimality) are still satisfied, the
following key fact is sufficient: 
\begin{itemize}
\item[] if $\{a,b\} \subseteq p$ or $\{a,b\} \subseteq q$ (in
  particular, if $(a,b)\in \dom(H_{\p \amal xy \q}) \cup 
\dom(F_{\p \amal  xy \q})$), 
\item[] then $\p \amal  xy \q \models a,b \< c$ holds iff at least one
  of the following is satisfied: 
\begin{enumerate}
\item $\p \models a,b \< c $
\item  $\q \models a,b \< c$
\item  $\p \models a,b \< x$, $\q \models y \< c$. 
\end{enumerate}
\end{itemize}
We leave the details of the argument to the reader. 
%
\end{proof}

\begin{Corollary}\label{plus}
Assume that $\p$, $\q$, $\r$ are as above, 
$\q\in \QQ_\delta$, $\p \le \M_\delta$.   

Then  $\p \oplus \q\in \QQ_\delta$ and 
 $\p \amal x y  \q \ \in \ \QQ_\delta$. 
\end{Corollary}

\begin{Definition}\label{hdef}
Let $(P,{\<})$ be a partial order of size $\kappa$, $\kappa$ 
a regular cardinal.    We say that  $P$  has \h0 if: 
\begin{quote}
Whenever $(\X_\alpha: \alpha <\kappa )$ is a sequence of finite 
suborderings of $(P,{\<})$,  
and $x_\alpha\in X_\alpha$ for all $\alpha < \kappa $, 
then there are $\alpha < {\beta}<\kappa $
such that:\\
$\X_\alpha$ and $\X_\beta$ agree on $X_\alpha\cap X_\beta$, 
$\X_\alpha$ and  $\X_\beta$ are separated 
end  extensions of $\X_\alpha \cap
\X_\beta$, $x_\alpha$ and $y_\alpha$ have the same type over $\X_\alpha
\cap \X_\beta$, and 
$${\X_\alpha} \amal {x_\alpha}{x_{\beta}}{\X_{\beta}}\le (P,{\<})$$  
\end{quote}
\end{Definition}

\begin{Remark}
By the $\Delta $-system lemma we know that for any sequence
$(X_\alpha: \alpha < \kappa)$ of finite sets there is a set $A
\subseteq \kappa$ of size $\kappa$ such that the sets $(X_\alpha:
\alpha \in A)$ form a $\Delta$-system.  If $(P, {\<})$ moreover
satisfies \ref{omega1like}, then we may additionally assume that for any
$\alpha < \beta < \kappa$ the models $\X_\alpha$ and $\X_\beta$ are 
separated end
extensions of $X_\alpha \cap \X_\beta$.  
So the important part of the above definition is really the last
clause. 
\end{Remark}

\begin{Lemma}\label{h0}
Assume that $(P,{\<})$ has  \h0. 
Then  $(P,{\<})$  satisfies conditions \ref{ccc} and  \ref{bounded}. 
\end{Lemma}

\begin{proof}
To show \ref{ccc}, consider any family 
$(x_\alpha : \alpha  < \kappa )$. 
 Let $\X_\alpha:= \{x_\alpha \}$. Since $(P,{\<})$ has
\h0, we can find $ \alpha < \beta$ such that 
$\X_\alpha \amal{x_\alpha }{x_\beta } \X_\beta  \le
(P,{\<})$, but this means that $ x_\alpha \<
x_\beta$. 

Finally we show \ref{bounded}. 
Let $f: (P,{\<}) \to (P,{\<})$ be monotone and assume that 
\ref{4prop}(\ref{bounded}) does not
hold.  This means that   for all small sets $A$ 
there is some $x$ such that $f(x) \not=x$
and $f(x)\notin A$.  So we can find a sequence
$(x_\alpha, y_\alpha: \alpha < \kappa)$ such that $\forall
\alpha\, x_\alpha \not= f(x_\alpha) = y_\alpha$ and 
the sets $\{x_\alpha, y_\alpha\}$ are pairwise disjoint.
  Let $\X_\alpha:= \{x_\alpha, y_\alpha \}\le (P,{\<})$. 
Wlog we either have $x_\alpha \< y_\alpha$ for all $\alpha$
or for no $\alpha$, similarly $y_\alpha \< x_\alpha$ for
all $\alpha$ or for no $\alpha$. 

Now find $\alpha$ and $\beta$ such that 
$\X_\alpha \amal{x_\alpha}{x_\beta} \X_\beta\le \M$. 
Now $\X_\alpha \amal{x_\alpha }{x_\beta } \X_\beta  \models
x_\alpha \< x_\beta$, but clearly 
 $\X_\alpha \amal{x_\alpha }{x_\beta }\X_\beta  \models
y_\alpha \not\< y_\beta$, so $f$ is not monotone, a
contradiction. 
\end{proof}

\section{Forcing}\label{forcing}

In this section we will carry out the forcing
construction that will prove the theorem.

We start with a universe $V_0$ satisfying GCH.  Let $\kappa =
{\omega_1} $,  and let $(M_{\delta}: {\delta} \le
\kappa)$,  $A_{\alpha {\beta}}^\delta $, \dots
${\Gamma}_{\alpha{\beta}}^\delta $ be as in \ref{setup}.

We define sequences $(\PP_\alpha,\name \QQ_\alpha: \alpha <
\kappa)$, 
$(\name \M_\alpha: \alpha < \kappa)$,
 $(\name R_\alpha: \alpha < {\kappa} )$ satisfying the
following for all $\alpha < \kappa$: 

\begin{enumerate}
\item $\PP_0 =\{\emptyset\}$, $\M_0 = \emptyset$.
\item $\name R_\alpha$ is a $\PP_\alpha$-name of a subset
        of $M_\alpha \times 
        M_\alpha $. 
\item  In $V^{\PP_\alpha}$, $F_{\alpha + 1} $ is constructed from
$R_\alpha$ as in \ref{construct}--\ref{cfact}.
\item $\name \QQ_\alpha$ is a $\PP_\alpha$-name, 
        $\forces_{\PP_\alpha }
         \name \QQ_\alpha = \QQ(\M_\alpha, F_{\alpha + 1} )$. 
\item $\PP_{\alpha+1} = \PP_\alpha * \name \QQ_\alpha$. 
\item $\PP_{\alpha +1} \forces \M_{\alpha + 1}$ is the
        \m0 defined by the generic filter on $\name
        \QQ_\alpha$, as in section \ref{code}. 
\item If $\alpha$ is a limit, then $\PP_\alpha$ is the
        finite support limit of $(\PP_\beta: \beta <
        \alpha)$, and 
        $\M_\alpha = \bigcup_{\beta < \alpha }
        \M_\beta$.  
\end{enumerate}
We let $\PP_{\kappa} $ be the finite support limit of
this iteration. 
Since all forcing notions involved satisfy  the countable chain
condition, we may (using the usual  bookkeeping arguments) assume
that: 
\begin{enumerate}
\item[$(*)$] Whenever $\name R$ is a
$\PP_{\kappa}$-name of a  subset of
$M_{\kappa} \times M_{\kappa}$ of size $< \kappa$, then there is some
$\alpha$ such that  
$\name R = \name R_\alpha$.  
\end{enumerate}

\begin{Lemma}\label{dense}
Let $A \subseteq M_{\delta}$ be finite.   Define a set 
$\bar \PP_{\delta}(A)$ as the set of all $\p\in \PP_{\delta}$ satisfying 
$$ \exists \M, \text{ $\M$ a finite \m0},
A \subseteq M, 
 \forall \alpha \in \dom(\p): 
\p\on \alpha \forces \p(\alpha)=\M\on M_{\alpha+1} $$
Then $\bar \PP_{\delta}(A)$ is dense in $\PP_{\delta}$.  In particular,
$\bar \PP_\delta := \bar \PP_\delta(\emptyset)$ is dense in $\PP_\alpha$. 
\end{Lemma}

\begin{proof}  By induction on ${\delta}$.   Limit steps are easy.  

Let $\p\in \PP_{\alpha+1}$.   By strengthening $\p\on \alpha$ we may
assume that there is a \m0 $\N$ such that $\p \on \alpha \forces
\p(\alpha)= \N$.  By \ref{dense1}(c) we may assume that
 $A \subseteq  N$.  By inductive assumption we 
may assume that there is a
\m0  $\M$ witnessing $\p \on \alpha \in \bar \PP_\alpha (N \cap
M_\delta)$.  Clearly $\M$ and $\N$ agree on $M \cap N$. 
  Define $\M':= \M \oplus  \N$.  So $\p\on \alpha \forces 
\p(\alpha) \weaker \M'$.  Define $\p' $ by
demanding $\p' \on \alpha =p\on \alpha$, $\p'(\alpha)=\M'$.  By
\ref{plus}, $\p'\on \alpha \forces \p'(\alpha)\in \QQ_\delta$. 
 Clearly
$\p' \stronger \p$, and $\p' \in \PP_{\alpha+1}$. 
\end{proof}

In $V^{\PP_{\kappa} }$, let $\M_{\kappa} = 
\bigcup_{\alpha < {\kappa} } \M_\alpha$.   We claim
that the structure $({\kappa}, \<_{\M_{\kappa} })$
satisfies the four conditions from \ref{4prop}.  The argument in
section \ref{code} shows that we have \ref{omega1like} and
\ref{defble}.  So, by \ref{h0} we only have to check that $\M$ has 
\h0. 

\medskip

Let 
$\name \X:= (\name \X_\alpha: \alpha < {\kappa})$ be a sequence of names for
finite \m0s, and let 
$(\name x_\alpha: \alpha < {\kappa})$ be forced to satisfy $\name x_\alpha
 \in X_\alpha$.  Assume that $\p$ is a condition forcing that
$\name\X$ witnesses the failure of \ref{hdef}. 
 For each $\alpha <{\kappa}$ we find $\p_\alpha \stronger \p $ which
decides $\X_\alpha$ 
and $x_\alpha$. Let $e_\alpha = \supp(\p_\alpha)$, $\delta_\alpha:=
\max(e_\alpha) + 1$, so $\p_\alpha \in \PP_{\delta_\alpha}$.

  By \ref{dense} we may
assume that there are finite \m0s $\P_\alpha$ such that $\forall
\varepsilon \in e_\alpha\,\, p_\alpha \on \varepsilon \forces
p(\varepsilon) = \P_\alpha \on {M_{\varepsilon{+}1}} $.

We may also assume that $X_\alpha \subseteq P_\alpha$, and hence,
(since $p_\alpha \forces \X_\alpha \le \M_\kappa$, 
$p_\alpha \forces \P_\alpha \le \M_\kappa$) $\X_\alpha \le \P_\alpha$. 

We may assume that the \m0s $\P_\alpha$ form a $\Delta$-system, say
with heart $\P^{\Delta} $ (so in particular $\P^{\Delta} \le
\P_\alpha$). Moreover, we may assume that each $\P_\alpha$ is a
separated  end
extension of $\P^\Delta $. 
Also, since the
$x_\alpha$ are wlog  all different, we may assume $x_\alpha \in
P_\alpha \setminus P^{\Delta}$.   Finally, we may assume that
each $x_\alpha$ has the same type over $\P^{\Delta}$.

Now pick any $\alpha < {\beta}$.  We will find a condition $\q
\stronger \p_\alpha$, $\q \stronger \p_{\beta}$ such that 
$\q \forces \P_\alpha \amal {x_\alpha}{x_{\beta}} \P_{\beta} \le
\M_\kappa $.

Let $\Q = \P_\alpha \amal {x_\alpha}{x_{\beta}} \P_{\beta} $. First
note that $\P_\alpha \le \Q$, $\P_{\beta} \le \Q$ by
\ref{amalfact}. Let $\varepsilon ^*$ be such that $x_\beta \in
M_{\varepsilon^*+1} \setminus M_\varepsilon$. 

Now note that 
$$\Q\on {M_{\varepsilon{+}1}} = 
\cases {\P_\alpha \on {M_{\varepsilon{+}1}} &
                        if $\varepsilon < \varepsilon^*$\cr
\P_\alpha \amal {x_\alpha}{x_{\beta}} \P_{\beta} \on
                {M_{\varepsilon{+}1}} & 
                        if $\varepsilon \ge \varepsilon^*$\cr}
$$

We now define a condition $\q$ with $\supp(\q) = 
\supp(\p_\alpha) \cup \supp(\p_{\beta})$, by stipulating 
$\q(\varepsilon) = \Q\on {M_{\varepsilon{+}1}} $. 

By \ref{plus} we know that $\q \on \varepsilon  \forces 
\Q\on M_{\varepsilon +1} \in \QQ_\varepsilon $, so by induction it is
clear that $\q$ is indeed a condition.  Clearly 
$\q \forces \P_\alpha \amal {x_\alpha}{x_{\beta}} \P_{\beta} \le
\M_\kappa $.



\end{document}